\documentclass[leqno]{article}

\usepackage{latexsym}
\usepackage{amssymb}
\usepackage{amsthm}
\usepackage{amsmath}

\usepackage{graphics}

\theoremstyle{definition}

\renewcommand{\L}{\ensuremath{\mathrm{L}}}

\newcommand{\D}{\ensuremath{{\cal D}}}
\renewcommand{\S}{\ensuremath{{\cal S}}}

\newcommand{\mb}[1]{\ensuremath{\mathbb{#1}}}
\newcommand{\N}{\mb{N}}

\newcommand{\R}{\mb{R}}

\newcommand{\cl}[1]{\ensuremath{\mathrm{cl}[#1]}}
\newcommand{\A}{\ensuremath{{\cal A}}}
\newcommand{\G}{\ensuremath{{\cal G}}}

\renewcommand{\d}{\ensuremath{\partial}}
\newcommand{\diff}[1]{\frac{d}{d#1}}

\newcommand{\grad}{\ensuremath{\mbox{\rm grad}\,}}

\newfont{\bl}{msbm10 scaled \magstep2}

\newtheorem{thm}{Theorem}
\newtheorem{lemma}[thm]{Lemma}
\newtheorem{prop}[thm]{Proposition}
\newtheorem{defn}[thm]{Definition}
\newtheorem{cor}[thm]{Corollary}
\newtheorem{rem}[thm]{Remark}
\newtheorem{ex}[thm]{Example}

\newcommand{\FT}[1]{\widehat{#1}}

\newcommand{\dis}[2]{\langle #1 , #2 \rangle}  
\newcommand{\inp}[2]{\langle #1 | #2 \rangle}  
\newcommand{\notmid}{\mid\kern-0.5em\not\kern0.5em}

\newcommand{\norm}[2]{{\| #1 \|}_{#2}}
\newcommand{\lone}[1]{\norm{#1}{\L^1}}

\newcommand{\linf}[1]{\norm{#1}{\L^\infty}}

\newcommand{\al}{\alpha}
\newcommand{\be}{\beta}
\newcommand{\ga}{\gamma}
\newcommand{\de}{\delta}
\newcommand{\eps}{\varepsilon}

\newcommand{\vphi}{\varphi}

\newcommand{\sig}{\sigma}

\newcommand{\supp}{\mathop{\mathrm{supp}}}

\newcommand{\ovl}[1]{\overline{#1}}

\newcommand{\Zin}{\ensuremath{C_*}}
\newcommand{\Zho}{\ensuremath{\dot{C}_*}}
\newcommand{\floor}[1]{\ensuremath{\lfloor #1 \rfloor}}
\newcommand{\inorm}[2]{\ensuremath{|#1|_{\Zin^{#2}}}}
\newcommand{\hnorm}[2]{\ensuremath{|#1|_{\Zho^{#2}}}}

\begin{document}

\title{Geophysical modelling with Colombeau functions: Microlocal properties
and Zygmund regularity}
\author{G\"{u}nther H\"{o}rmann and Maarten V. de Hoop \\
        \textit{
        Department of Mathematical and Computer Sciences}, \\
        \textit{
        Colorado School of Mines, Golden CO 80401}}
\date{\today}
\maketitle

\begin{abstract}
  In global seismology Earth's properties of fractal nature occur.
  Zygmund classes appear as the most appropriate and systematic way to
  measure this local fractality. For the purpose of seismic wave
  propagation, we model the Earth's properties as Colombeau
  generalized functions. In one spatial dimension, we have a precise
  characterization of Zygmund regularity in Colombeau algebras. This
  is made possible via a relation between mollifiers and wavelets.
\end{abstract}

\section{Introduction}

\emph{Wave propagation in highly irregular media}. In global
seismology, (hyperbolic) partial differential equations the
coefficients of which have to be considered generalized functions; in
addition, the source mechanisms in such application are highly
singular in nature. The coefficients model the (elastic) properties of
the Earth, and their singularity structure arises from geological and
physical processes. These processes are believed to reflect themselves
in a multi-fractal behavior of the Earth's properties. Zygmund classes
appear as the most appropriate and systematic way to measure this
local fractality (cf.~\cite[Chap.4]{Holschneider:95}).

\emph{The modelling process and Colombeau algebras}. In the seismic
transmission problem, the diagonalization of the first order system of
partial differential equations and the transformation to the second order
wave equation requires differentiation of the coefficients. Therefore,
highly discontinuous coefficients will appear naturally although the
original model medium varies continuously.  However, embedding the fractal
coefficient first into the Colombeau algebra ensures the equivalence after
transformation and yields unique solvability if the regularization scaling
$\ga$ is chosen appropriately (cf.~\cite{LO:91,O:89,HdH:01}). We use the
framework and notation (in particular, $\G$ for the algebra and $\A_N$ for
the mollifier sets) of Colombeau algebras as presented in \cite{O:92}.

An interesting aspect of the use of Colombeau theory in wave propagation is
that it leads to a natural control over and understanding of `scale'. In
this paper, we focus on this modelling process.

\section{Basic definitions and constructions}

\subsection{Review of Zygmund spaces}

We briefly review homogeneous and inhomogeneous
Zygmund spaces, $\Zho^s(\R^m)$ and $\Zin^s(\R^m)$, via a characterization
in pseudodifferential operator style which follows essentially the
presentation in \cite{Hoermander:97}, Sect.~8.6. Alternatively, for practical
and implementation issues one may prefer the characterization via growth
properties of the discrete wavelet transform using orthonormal wavelets
(cf.~\cite{Meyer:92}).

Classically, the Zygmund spaces were defined  as extension of H\"older spaces
by boundedness properties of difference quotients. Within the systematic and
unified approach of Triebel (cf.~\cite{Triebel:I,Triebel:II}) we can simply
identify the Zygmund spaces in a scale of inhomogeneous and homogeneous
(quasi) Banach spaces, $B^s_{p q}$ and $\dot{B}^s_{p q}$ ($s\in\R$,
$0 < p, q \leq \infty$), by $\Zin^s(\R^m) = B^s_{\infty \infty}(\R^m)$
and $\Zho^s(\R^m) = \dot{B}^s_{\infty \infty}$. Both $\Zin^s(\R^m)$ and
$\Zho^s(\R^m)$ are Banach spaces.

To emphasize the close relation with mollifiers we describe a characterization
of Zygmund spaces in pseudodifferential operator style in more detail.

Let $0 < a < b$ and choose $\vphi_0\in\D(\R)$, $\vphi_0$ symmetric and
positive, $\vphi_0(t) = 1$ if $|t| < a$, $\vphi_0(t) = 0$ if $|t| > b$, and
$\vphi_0$ strictly decreasing in the interval $(a,b)$. Putting
$\vphi(\xi) = \vphi_0(|\xi|)$ for $\xi\in\R^m$ then defines a function
$\vphi\in\D(\R^m)$. Finally we set
\[
        \psi(\xi) = - \inp{\xi}{\grad\vphi(\xi)}
\]
and note that if $a < |\xi| < b$ then
$\psi(\xi) = - \vphi_0'(|\xi|) |\xi| > 0$. We denote by ${\cal M}(\R^m)$ the
set of all pairs $(\vphi,\psi)\in\D(\R^m)^2$
that are constructed as above (we usually suppress the dependence of ${\cal
M}$ on $a$ and $b$ in the notation).

We are now in aposition to state the characterization theorem for the
inhomogeneous Zygmund spaces as subspaces of $\S'(\R^m)$. It follows from
\cite{Triebel:88}, Sec.~2.3, Thm.~3 or, alternatively, from
\cite{Hoermander:97}, Sec.~8.6. Note that all appearing pseudodifferential
operators in the following have $x$-independent symbols and are thus given
simply by convolutions.

\begin{thm}\label{inh_Z} Assume that $a \leq 1/4$ and $b \geq 4$ and choose
$(\vphi,\psi)\in {\cal M}(\R^m)$ arbitrary. Let $s\in\R$ then
$u\in\S'(\R^m)$ belongs to the inhomogeneous
Zygmund space of order $s$ $\Zin^s(\R^m)$ if and only if
\begin{equation}
  \inorm{u}{s} := \linf{\vphi(D)u} +
                \sup\limits_{0 < t < 1}\Big( t^{-s} \linf{\psi(tD)u}\Big) <
         \infty   .
\end{equation}
\end{thm}
(Note that we made use of the modification for $q=\infty$ in
\cite{Triebel:88}, equ.~(82).)

\begin{rem}\label{Z_rem} 
\begin{enumerate}
\item $\inorm{u}{s}$ defines an equivalent norm on $\Zin^s$. In fact that
        all norms defined as above by some $(\vphi,\psi)\in{\cal M}(\R^m)$ are
        equivalent can be seen as in \cite{Hoermander:97}, Lemma 8.6.5.
\item If $s\in \R_+ \setminus \N$ then $C^s_*(\R^m)$ is the classical
  H\"older space of regularity $s$. Denoting by $\floor{s}$ the greatest
  integer less than $s$ it consists of all $\floor{s}$ times  continuously
  differentiable functions $f$ such that $\d^\al f$ is bounded when $|\al|
  \leq \floor{s}$ and globally H\"older continuous with exponent
  $s-\floor{s}$ if $|\al| = \floor{s}$.
\item Due to the term $\linf{\vphi(D)u}$ the norm $\inorm{u}{s}$ is not
        homogeneous with respect to a scale change in the argument of $u$.
\item If $u\in\L^\infty(\R^m)$ then (cf.~\cite{Hoermander:97}, Sect.~8.6)
\begin{equation}
        u(x) = \vphi(D)u(x) + \int\limits_1^\infty \psi(D/t)u(x) \frac{dt}{t}
        \qquad \text{ for almost all } x .
\end{equation}
Using $\vphi(\xi) = \int_0^1 \psi(\xi/t)/t \,dt$ this can be rewritten in
the form $u(x) = \int_0^\infty \psi(D/t)u(x)/t \,dt$
and
resembles Calderon's classical identity in terms of a continuous wavelet
transform (cf.~\cite{Meyer:92}, Ch.~1, (5.9) and (5.10)).
\item In a similar way one can characterize the homogeneous Zygmund spaces as
subspaces of $\S'(\R^m)$ modulo the polynomials ${\cal P}$. A proof
can be found in \cite{Triebel:82}, Sec.~3.1, Thm.~1. We may identify
$\S'/{\cal P}$ with the dual space $\S_0'(\R^m)$ of
$\S_0(\R^m) = \{ f\in\S(\R^m) \mid \d^\al \FT{f}(0) = 0 \,
\forall \al\in\N_0^m \}$, the Schwartz functions with
vanishing moments, by mapping the class $u+{\cal P}$ with representative
$u\in\S'$ to $u\mid_{\S_0}$. Assume that $a \leq 1/4$ and $b \geq 4$ and
choose $\psi\in\D(\R^m)$ as constructed above and let $s\in\R$ and
$u\in\S'(\R^m)$. Then $u\!\!\mid_{\S_0}$ belongs to the
homogeneous Zygmund space $\Zho^s(\R^m)$ of order $s$  if and only if
\begin{equation}
  \hnorm{u}{s} :=
                \sup\limits_{0 < t < \infty}\Big( t^{-s} \linf{\psi(tD)u}\Big) <
         \infty   .
\end{equation}
(Note that we use the modification for $q=\infty$ in \cite{Triebel:82},
equ.~(16).)
\end{enumerate}
\end{rem}

\subsection{The continuous wavelet transform}

Following \cite{Holschneider:95} we call a function
$g\in\L^1(\R^m)\cap\L^\infty(\R^m)$ with $\int g = 0$ a \emph{wavelet}. We
shall say that it is a \emph{wavelet of order $k$} ($k\in\N_0$) if the
moments up to order $k$ vanish, i.e., $\int x^\al g(x) dx = 0$ for
$|\al|\leq k$.

The (continuous) wavelet transform is defined for $f\in\L^p(\R^m)$
($1\leq p \leq \infty$) by ($\eps > 0$)
\begin{equation}\label{wf_trafo}
  W_g f(x,\eps) = \int\limits_{\R^m} f(y) \frac{1}{\eps^m}
  \bar{g}(\frac{y-x}{\eps}) \, dy = f * (\bar{\check{g}})_\eps (x)
\end{equation}
where we have used the notation $\check{g}(y) = g(-y)$ and $g_\eps(y) =
g(y/\eps)/\eps^n$. By Young's inequality $W_g f(.,\eps)$ is in $\L^p$ for all
$\eps > 0$ and $W_g$ defines a continuous operator on this space for each
$\eps$.

If $g\in C_c(\R^m)$ we can define $W_g f$ for $f\in\L^1_{\text{loc}}(\R^m)$
directly by the same formula (\ref{wf_trafo}). If $g\in\S_0(\R^m)$ then $W_g$
can be extended to $\S'(\R^m)$ as the adjoint of the wavelet synthesis (cf.\
\cite{Holschneider:95}, Ch.~1, Sects.\ 24, 25, and 30) or directly by
$\S'$-$\S$-convolution in formula (\ref{wf_trafo}).

\begin{rem}
If $f$ is a polynomial and $g\in\S_0$ it is easy to see that $W_g f = 0$.
In fact,  $f$, $g$, and $W_g f$ are in $\S'$ and $\FT{(W_g f(.,\eps))} =
\FT{f} \bar{\FT{g}}(\eps .)$. Since $g$ is in $\S_0$ the Fourier transform
$\FT{g}(\eps .)$ is smooth and vanishes of infinite order at $0$. But
$\FT{f}$ has to be a linear combination of derivatives of $\de_0$ implying
$\FT{f} \bar{\FT{g}}(\eps .) = 0$. Therefore the wavelet transform `is
blind to polynomial parts' of the analyzed function (or distribution) $f$.
In terms of geophysical modelling this means that a polynomially varying
background medium is filtered out automatically.
\end{rem}

\subsection{Wavelets from mollifiers}

The Zygmund class characterization in Theorem \ref{inh_Z} (and remark
\ref{Z_rem},(v)) used asymptotic estimates of scaled smoothings of the distribution which
resembles typical mollifier constructions in Colombeau theory. In this
subsection we relate this in turn directly to the wavelet transform
obtaining the well-known wavelet characterization of Zygmund spaces.

Let $\chi\in\S(\R^m)$ with $\int \chi = 1$ and define the function $\mu$ by
\begin{equation}\label{mo_to_wv}
 \ovl{\check{\mu}(x)} := m \chi(x) + \inp{x}{\grad\chi(x)} \, .
\end{equation}
Then $\mu$ is in $\S(\R^m)$ and is a wavelet since a simple integration by
parts shows that
\begin{multline*}
  (-1)^{|\al|} \ovl{\int \mu(x) x^\al \, dx} =
    \int \bar{\check{\mu}}(x) (-x)^\al \, dx \\
= (-1)^{|\al|+1} \int x^\al \chi(x)\, dx  \sum_{j=1}^m \al_j
 =  (-1)^{|\al|+1} |\al| \int x^\al \chi(x)\, dx \; .
\end{multline*}
$\int \mu = 0$ and if $|\al| > 0$ we have $\int x^\al \mu(x)\,dx = 0$
if and only if $\int x^\al \chi(x) \,dx = 0$.
Therefore $\mu$ defined by (\ref{mo_to_wv}) is a wavelet of order $N$ if and
only if the mollifier $\chi$ has vanishing moments of order $1 \leq |\al|
\leq N$.

Furthermore, by straightforward computation, we have
\begin{equation}\label{mo_to_wv_2}
 (\bar{\check{\mu}})_\eps(x) =
  -\eps \diff{\eps} \big(\chi_\eps (x)\big)
\end{equation}
yielding an alternative of (\ref{mo_to_wv}) in the form
$\bar{\check{\mu}}(x) = - \diff{\eps}\big(\chi_\eps(x)\big)\mid_{\eps=1}$.

If $(\vphi,\psi)\in{\cal M}(\R^m)$ arbitrary and $\chi$, $\mu$ are the unique Schwartz
functions such that $\FT{\chi} = \vphi$ and $\FT{\mu} = \psi$, then
straightforward computation shows that $\mu$ and $\chi$ satisfy the relation
(\ref{mo_to_wv}). Therefore since $\mu$ is then a real valued and even
wavelet we have for $u\in\S'$
\[
        \psi(tD)u(x) = t^m \bar{\check{\mu}}(\frac{.}{t})*u(x) = W_\mu u(x,t) \;
        .
\]
Hence the distributions $u$ in the Zygmund class $\Zin^s(\R^m)$ can be
characterized in terms of a wavelet transform and a smoothing
pseudodifferential operator by
$ \linf{\vphi(D)u} < \infty$  and
$ \sup_{0 < t < 1} \Big( t^{-s} \linf{W_\mu u(.,t)}\Big) < \infty$.
We have shown
\begin{thm} Let $(\FT{\chi},\FT{\mu})\in{\cal M}(\R^m)$. A distribution
$u\in\S'(\R^m)$ belongs to the Zygmund class $\Zin^s(\R^m)$ if and only if
\begin{equation}\label{Z_W_char}
  \linf{u * \chi} < \infty \quad \text{and} \quad
                \linf{W_\mu u(.,r)} = O(r^s) \;\; (r \to 0) .
\end{equation}
\end{thm}

\begin{rem}\label{W_rem}
\begin{enumerate}
\item Observe that the condition on $\FT{\chi}$ implies that $\chi$ and hence
$\mu$ can never have compact support. If this characterization is to be
used in a theory of Zygmund regularity detection within Colombeau algebras
one has to allow for mollifiers of this kind in the corresponding embedding
procedures. This is the issue of the following subsection. Nevertheless we
note here that according to remarks in \cite[(2.2) and (3.1)]{Jaffard:97a}
and, more precisely, in \cite[Ch.3]{Meyer:98} the restrictions on the
wavelet itself in a characterization of type (\ref{Z_W_char}) may be
considerably relaxed
--- depending on the generality one wishes to allow for the analyzed
distribution $u$. However, in case $m=1$ and $u$ a function a flexible and
direct characterization (due to Holschneider and Tchamitchian) can be found
in \cite{Daubechies:92}, Sect.~2.9, or \cite{Holschneider:95}, Sect.~4.2.
\item There are more refined results in the spirit of the above theorem
        describing local H\"older (Zygmund) regularity by growth properties of
        the  wavelet transform (cf.~in particular \cite{Holschneider:95},
        Sect.~4.2, \cite{JM:96}, and \cite{Jaffard:97a}).
\item The counterpart of (\ref{Z_W_char}) for $\L^1_{\text{loc}}$-functions
        in terms of (discrete) multiresolution approximations is
        \cite{Meyer:92}, Sect.~6.4, Thm.~5.
\end{enumerate}
\end{rem}

\section{Colombeau modelling and wavelet transform}

\subsection{Embedding of temperate distributions}

We consider a variant of the Colombeau embedding
$\iota_\chi^\ga : \D'(\R^m) \to \G(\R^m)$ that was discussed in
\cite{HdH:01}, subsect.~3.2. As indicated in remark \ref{W_rem},(i)
we need to allow for mollifiers with noncompact support in order to gain
the flexibility of using wavelet-type arguments for the extraction of
regularity properties from asymptotic estimates. On the side of the embedded
distributions this forces us to restrict to $\S'$, a space still large
enough for the geophysically motivated coefficients in model PDEs.

Recall (\cite{HdH:01}, Def.~11) that an admissible scaling is
defined to be a continuous function $\ga : (0,1) \to \R_+$ such that
$\ga(r) = O(1/r)$, $\ga(r) \to\infty$, and $\ga(sr) = O(\ga(r))$ if $0<s<1$
(fixed) as $r\to 0$.
\begin{defn}
Let $\ga$ be an admissible scaling, $\chi\in\S(\R^m)$ with
$\int \chi = 1$, then we define
$\iota_\chi^\ga : \S'(\R^m) \to \G(\R^m)$ by
\begin{equation}
  \iota_\chi^\ga(u) = \cl{(u * \chi^\ga(\phi,.))_{\phi\in\A_0(\R^m)}} \qquad
   u\in\S'(\R^m)
\end{equation}
where
\begin{equation}
        \chi^\ga(\phi,x) = \ga(l(\phi_0))^m \chi(\ga(l(\phi_0)) x)
        \quad \text{ if } \phi = \phi_0\otimes \cdots \otimes\phi_0
        \quad  \text{ with } \phi_0\in\A_0 .
\end{equation}
\end{defn}
$\iota_\chi^\ga$ is well-defined since $(\phi,x) \to u*\chi^\ga(\phi,x)$
is clearly moderate and negligibility is preserved under this scaled
convolution. By abuse of notation we will write $\iota_\chi^\ga(u)(\phi,x)$
for the standard representative of $\iota_\chi(u)$.

The following statements describe properties of such a modelling procedure
resembling the original properties used by M. Oberguggenberger in
\cite{O:89}, Prop.1.5, to ensure unique solvability of symmetric hyperbolic
systems of PDEs (cf.~\cite{O:89,LO:91}). The definition of Colombeau
functions of logarithmic and bounded type is given in \cite{O:92},
Def.~19.2, the variation used below is an obvious extension.
\begin{prop} \leavevmode
\begin{enumerate}
\item $\iota_\chi^\ga : \S'(\R^m) \to \G(\R^m)$ is linear, injective, and
        commutes with partial derivatives.
\item $\forall u\in\S'(\R^m)$: $\iota_\chi^\ga(u) \approx u$.
\item If $u\in W^{-1}_\infty(\R^m)$ then $\iota_\chi^\ga(u)$ is of
        \emph{$\ga$-type}, i.e., there is $N\in\N_0$ such that for all
        $\phi\in\A_N(\R^m)$ there exist $C > 0$ and $1 > \eta > 0$:
\begin{equation}
        \sup\limits_{y\in\R^m} | \iota_\chi^\ga(u)(\phi_\eps,y)| \leq
                N \ga(C\eps) \quad 0 < \eps < \eta .
\end{equation}
\item If $u\in\L^\infty(\R^m)$ then $\iota_\chi^\ga(u)$ is
        of bounded type and its first order derivatives are of $\ga$-type.
\end{enumerate}
\end{prop}
\begin{proof}
\emph{ad (i),(ii):} Is clear from $\chi_\eps :=
        \chi^\ga(\phi_\eps,.) \to \de$ in $\S'$ as $\eps\to 0$ and the
        convolution formula.

\emph{ad (iii):} Although this involves only marginal changes in the proof of
        \cite{O:89}, Prop.~1.5(i), we recall it here to make the
        presentation more self-contained.

        Let $u = u_0 + \sum_{j=1}^m \d_j u_j$ with $u_j\in\L^\infty$
        ($j=0,\ldots,m$)        then with $\ga_\eps := \ga(\eps l(\phi_0))$
\begin{multline*}
        |u*\chi_\eps(x)| \leq \linf{u_0 * \chi_\eps} +
                \sum_{j=1}^m \linf{u_j * \d_j(\chi_\eps)} \\
        \leq \linf{u_0} \lone{\chi} +
                \ga_\eps \sum_{j=1}^m \linf{u_j} \lone{\d_j\chi} \\
                = \ga_\eps \big( \frac{\linf{u_0} \lone{\chi}}{\ga_\eps} +
                                \sum_{j=1}^m \linf{u_j} \lone{\d_j \chi} \big)
\end{multline*}
where the expression within brackets on the r.h.s.\ is bounded by some
constant $M$, dependent on $u$ and $\chi$ only but independent of $\phi$,
as soon as $\eps < \eta$ with $\eta$ chosen appropriately (and dependent
on $M$, $u$, $\chi$, and $\phi$). Therefore the assertion is
proved by putting $N \geq M$ and $C = l(\phi_0)$.

\emph{ad (iv):} is proved by similar reasoning
\end{proof}

In particular, we can model a fairly large class of distributions as
Colombeau functions of logarithmic growth (or log-type) thereby ensuring
unique solvability of hyperbolic PDEs incorporating such as coefficients.
\begin{cor} \leavevmode
\begin{enumerate}
\item If $\ga(\eps) = \log(1/\eps)$ then
        $\iota_\chi^\ga(W^{-1,\infty}) \subseteq \{ U\in\G \mid U
        \text{ is of log-type } \}$
        and
        \[
                \iota_\chi^\ga(\L^\infty) \subseteq
                \{ U\in\G \mid U \text{ of bounded type and }
                  \d^\alpha U \text{ of log-type for } |\al| = 1 \} .
        \]
\item If $u\in W^{-k,\infty}(\R^m)$ for $k\in\N_0$ then $\iota_\chi^\ga(u)$
        is of $\ga^k$-type. In particular, there is an admissible scaling $\ga$
        such that $\iota_\chi^\ga(u)$ and all first order derivatives
        $\d_j \iota_\chi^\ga(u)$ ($j = 1,\ldots,m$) are of log-type.
\end{enumerate}
\end{cor}

\subsection{Wave front sets under the embedding}

One of the most important properties of the embedding procedure introduced
in \cite{HdH:01} was its faithfulness with respect to the microlocal
properties if `appropriately measured' in terms of the set of $\ga$-regular
Colombeau functions $\G_\ga^\infty(\R^m)$ (\cite{HdH:01}, Def.~11). But
there the proof of this microlocal invariance property heavily used the
compact support property of the standard mollifier $\chi$ which is no longer
true in the current situation. In this subsection we show  how to extend the
invariance result to the new embedding procedure defined above.

\begin{thm} Let $w\in\S'(\R^m)$, $\ga$ an admissible scaling, and
$\chi\in\S(\R^m)$ with $\int \chi = 1$ then
\begin{equation}
        WF_g^\ga(\iota_\chi^\ga(w)) = WF(w) .
\end{equation}
\end{thm}
\begin{proof}
The necessary changes in the proof of \cite{HdH:01}, Thm.~15, are
minimal once we established the following
\begin{lemma} If $\vphi\in\D(\R^m)$ and $v\in\S'(\R^m)$ with $\supp(\vphi)
        \cap \supp(v) = \emptyset$ then
        $\vphi\cdot\iota_\chi^\ga(v)\in\G_\ga^\infty$.
\end{lemma}
\begin{proof} Using the short-hand notation
$\chi_\eps = \chi^\ga(\phi_\eps,.)$ and $\ga_\eps = \ga(\eps l(\phi_0))$ we
have
\[
        \d^\be\big(\vphi (v*\chi_\eps)\big)(x) =
                \ga_\eps^m \sum_{\al\leq\be} \binom{\be}{\al} \d^{\be-\al}\vphi(x) \,
                 \ga_\eps^{|\al|}\,     \dis{v}{\d^{\al}\chi(\ga_\eps(x-.))} \; .
\]
Hence we need to estimate terms of the form
$\ga_\eps^{|\al|}\, \dis{v}{\d^{\al}\chi(\ga_\eps(x-.))}$ when
$x\in\supp(\vphi) =: K$. Let $S$ be a closed set satisfying $\supp(v)
\subset S \subset \R^m \setminus K$ and put $d = \text{dist}(S,K)>0$.
Since $v$ is a temperate distribution there is $N\in\N$ and $C > 0$ such
that
\[
  \ga_\eps^{|\al|} |\dis{v}{\d^{\al}\chi(\ga_\eps(x-.))}| \leq
        C  \ga_\eps^{|\al|} \sum_{|\sig|\leq N} \sup\limits_{y\in S}
                |\d^\sig\big( \d^\al \chi(\ga_\eps(x-y)) \big)|  \, .
\]
$\chi\in\S$ implies that each term in the sum on the right-hand side can be
estimated for arbitrary $k\in\N$  by
\[
  \sup\limits_{y\in S} |\d^{\sig+\al}\chi(\ga_\eps(x-y)) \big)|
          \ga_\eps^{|\sig|} \leq
        \ga_\eps^{|\sig|} \sup\limits_{y\in S} C_k (1+\ga_\eps|x-y|)^{-k}
        \leq C'_k \ga_\eps^{|\sig|-k}/ d^k
\]
if $x$ varies in $K$. Since $|\al|+|\sig| \leq |\be|+ N$ we obtain
\[
        \linf{\d^\be\big(\vphi (v*\chi_\eps)\big)} \leq
                C' \ga_\eps^{m+N +|\be|-k}
\]
with a constant $C'$ depending on $k$, $v$, $\vphi$, $d$, and $\chi$ but
$k$ still arbitrary. Choosing $k = |\be|$, for example, we conclude that
$\vphi\cdot\iota_\chi^\ga(v)$ has a uniform $\ga_\eps$-growth over all
orders of derivatives. Hence it is a $\ga_\eps$-regular Colombeau function.
\end{proof}

Referring to the proof (and the notation) of \cite{HdH:01}, Thm.~15, we may
now finish the proof of the theorem simply by carrying out the following
slight changes in the two steps of that proof.

\emph{Ad step 1:} Choose $\psi\in\D$ such that $\psi = 1$ in a neighborhood
of $\supp(\vphi)$ and write
\[
        \vphi (w*\chi_\eps) = \vphi ((\psi w)*\chi_\eps) + \vphi
                (((1-\psi)w)*\chi_\eps \; .
\]
The first term on the right can be estimated by the same methods as
in \cite{HdH:01} and the second term is $\ga$-regular by the lemma above.

\emph{Ad step 2:} Rewrite
\[
        \vphi w = \vphi \psi w = \vphi (\psi w - (\psi w)*\chi_\eps) +
                \vphi((\psi w)*\chi_\eps)
\]
and observe that the reasoning of \cite{HdH:01} is applicable since
$\Sigma_g^\ga(\vphi\, \iota_\chi^\ga(\psi w)) \subseteq
\Sigma_g^\ga(\vphi\, \iota_\chi^\ga(w))$ by the above lemma.

\end{proof}

\subsection{The modelling procedure and wavelet transforms}

Simple wavelet-mollifier correspondences as in subsection 2.3 allow us to
rewrite the Colombeau modelling procedure and hence prepare for the
detection of original Zygmund regularity in terms of growth properties in
the scaling parameters.

A first version describes directly $\iota_\chi^\ga$ but involves an
additional nonhomogeneous term.
\begin{lemma}\label{inhom_lemma}
If $\chi\in\S(\R^m)$ has the properties $\int\chi = 1$ and
$\int x^\al \chi(x) dx = 0$  ($0 < |\al| \leq N$) then $\bar{\check{\mu}} =
 -\diff{\eps}(\chi_\eps)\mid_{\eps=1}$ defines a wavelet of order $N$
 and  we have for any $f\in\S'(\R^m)$
\begin{equation}
 \iota_\chi^\ga(f)(\phi,x) = f * \chi (x) + \!\!\!\!
  \int\limits_{1/\ga(l(\phi))}^{1}\!\!\!\! W_\mu f(x,r)\, \frac{dr}{r} \;.
\end{equation}
\end{lemma}
\begin{proof}
Let $\eps > 0$ then eq.~(\ref{mo_to_wv_2}) implies  $W_\mu f(x,\eps) =
f * (\bar{\check{\mu}})_\eps (x) =
-\eps \diff{\eps}\big( f * \chi_\eps(x) \big)$
and integration with respect to $\eps$ from $1/\ga(l(\phi))$ to $1$ yields
\[
  - \!\!\!\! \int\limits_{1/\ga(l(\phi))}^{1}
  \!\!\!\! W_\mu f(x,\eps) \frac{d\eps}{\eps} =  f * \chi (x) -
  \iota_\chi^\ga(f)(\phi,x) \; .
\]
\end{proof}

A more direct mollifier wavelet correspondence is possible via
derivatives of $\iota_\chi^\ga$ instead.
\begin{lemma}\label{hom_lemma}
If $\chi\in\S(\R^m)$ with $\int \chi =1$ then for any $\al\in\N_0^n$ with
$|\al| > 0$
\begin{equation}\label{chi_al}
  \chi_\al(x) = \ovl{(\d^\al\chi)\check{\ }(x)}
\end{equation}
is a wavelet of order $|\al|-1$ and for any $f\in\S'(\R^m)$ we have
\begin{equation}
 \d^\al \iota_\chi^\ga(f)(\phi,x) =
   \ga(l(\phi))^{|\al|}\, W_{\chi_\al} f (x,\frac{1}{\ga(l(\phi))}) \; .
\end{equation}
\end{lemma}
\begin{proof}
Let $|\be| < |\al|$ then $\int x^\be D^\al\chi(x) \, dx =
(-D)^\be(\xi^\al \FT{\chi}(\xi))\mid_{\xi=0} = 0$
which proves the first assertion. The second assertion follows from
\[
 \d^\al \iota_\chi^\ga(f)(\phi) = \ga^{|\al|+m} f *
 \d^\al\chi(\ga .) =
  \ga^{|\al|} f *
  \ovl{\big(\ga^m (\d^\al\bar{\chi})\check{\ }(\ga .)\big)\check{\ }}
\]
with the short-hand notation $\ga = \ga(l(\phi))$.
\end{proof}

Both lemmas \ref{inhom_lemma} and \ref{hom_lemma} may be used to translate
(global) Zygmund regularity of the modeled (embedded) distribution $f$ via
Thm.~\ref{Z_W_char} into asymptotic growth properties with respect to the
regularization parameter. To what extent this can be utilized to develop a
faithful and completely intrinsic Zygmund regularity theory of Colombeau
functions may be subject of future research.

\section{Zygmund regularity of Colombeau functions: the one-dimensional
         case}

If we combine the basic ideas of the Zygmund class characterization in 2.3
with the simple observations in 3.3 we are naturally lead to define a
corresponding regularity notion intrinsically in Colombeau algebras as
follows.
\begin{defn}\label{Z_C_def} Let $\ga$ be an admissible scaling function and
$s$ be a real number.
A Colombeau function $U\in\G(\R^m)$ is said to be \emph{globally of
$\ga$-Zygmund regularity $s$} if for all $\alpha\in\N_0^m$ there is $M\in\N_0$
such that for all $\phi\in\A_M(\R^m)$ we can find positive constants $C$ and
$\eta$ such that
\begin{equation}\label{ZC_def}
        |\d^\al U(\phi_\eps,x)| \leq
            \begin{cases}
               C & \text{if $|\al| < s$}\\
                C \;\ga_\eps^{|\al|-s} & \text{if $|\al| \geq s$}
            \end{cases}
            \qquad
        x\in\R^m, 0<\eps<\eta .
\end{equation}
The set of all (globally) $\ga$-Zygmund regular Colombeau functions of
order $s$ will be denoted by $\G_{*,\ga}^s(\R^m)$.
\end{defn}

A detailed analysis of $\G_{*,\ga}^s$ in arbitrary space dimensions and not
necessarily positive regularity $s$ will appear elsewhere. Here, as an
illustration, we briefly study the case $m=1$ and $s>0$ in some detail.
Concerning applications to PDEs this would mean that we are allowing for
media of typical fractal nature varying continuously in one space
dimension. For example one may think of a coefficient function $f$ in
$\Zin^s(\R)$ to appear in the following ways.
\begin{ex}
\begin{enumerate}
\item Let $f$ be constant outside some interval $(-K,K)$ and
        equal to a typical trajectory of Brownian motion in $[-K,K]$; it is
        well-known that with probability $1$ those trajectories are in $\Zho^s(\R)$
        whenever $s < 1/2$. This is proved, e.g., in \cite{Holschneider:95}, Sect.~4.4,
        elegantly by wavelet transform methods.
\item We refer to \cite{Zygmund:68}, Sect.~V.3, for notions and notation in
        this example. Then similarly to the above one can set $f=0$ in
        $(-\infty,0]$, $f=1$ in $[2\pi,\infty)$ and in $[0,2\pi]$ let $f$ be
        Lebesgue's singular function associated with a Cantor-type set of order
        $d\in\N$ with (constant) dissection ratio $0 < \xi < 1/2$. Then $f$
        belongs to $\Zin^s(\R)$ with $s = \log(d+1)/|\log(\xi)|$. (The classical
        triadic Cantor set corresponds to the case $d=2$ and $\xi = 1/3$.)
\end{enumerate}
\end{ex}

We have already seen that the Colombeau embedding does not change the
microlocal structure (i.e., the $\ga$-wave front set) of the original
distribution. We will show now that also the refined Zygmund regularity
information is accurately preserved. If $n\in\N_0$ we denote by
$C^n_b(\R)$ the set of all $n$ times continuously differentiable functions with
the derivatives up to order $n$ bounded. Note that $C^n_b(\R)$ is a strict
superset of $\Zin^{n+1}(\R)$.
\begin{thm} Let $\chi\in\A_0(\R)$ and $s>0$. Define $n\in\N_0$ such that
$n < s \leq n+1$ then we have
\begin{equation}
        \iota_\chi^\ga(C^n_b(\R)) \cap \G_{*,\ga}^s(\R) =
        \iota_\chi^\ga(\Zin^s(\R)) .
\end{equation}
In other words, in case $0< s <1$ we can  precisely identify those Colombeau
functions that arise from the Zygmund class of order $s$ within all
embedded bounded continuous functions.
\end{thm}
\begin{proof} We use the characterizations in \cite{Daubechies:92}, Thms.~2.9.1
and 2.9.2 and the remarks on p.~48 following those; choosing a smooth
compactly supported wavelet $g$ of order $n$ we may therefore state the
following\footnote{Note that we do not use the
wavelet scaling convention adapted to $\L^2$-spaces here}: $f\in C^n_b(\R)$
belongs to $\Zin^s(\R)$ if and only if there is $C>0$ such that
\begin{equation}\label{Daub}
        |W_g f(x,r)| \leq C r^s \quad \text{ for all } x .
\end{equation}

Now the proof is straightforward. First let $f\in\Zin^s(\R)$. If
$|\al| - 1 < n$ then $|\d^\al \iota_\chi^\ga(f)| =
|\iota_\chi^\ga(\d^\al f)| \leq \linf{\d^\al f} \lone{\chi}$ by Young's
inequality. If $|\al| > n$ we use lemma \ref{hom_lemma} and set $\ga_\eps =
\ga(\eps l(\phi))$ to obtain
\[
        \d^\al \iota_\chi^\ga(f)(\phi_\eps,x) = \ga_\eps^{|\al|} W_{\chi_\al}
        f(x,\ga_\eps^{-1})
\]
where $\chi_\al$ is a wavelet of order at least $n$. Hence (\ref{Daub})
gives an upper bound $C \ga_\eps^{|\al| - s}$ uniformly in $x$. Hence
(\ref{ZC_def}) follows.

Finally, if we know that $f\in C^n_b(\R)$ and
$\iota_\chi^\ga(f) \in \G_{*,\ga}^s(\R)$ then combination of (\ref{ZC_def})
and lemma \ref{hom_lemma} gives if $|\al| \geq s$
\[
        |\ga_\eps^{|\al|} W_{\chi_\al} f(x,\ga_\eps^{-1})| =
                |\d^\al \iota_\chi^\ga(f)(\phi_\eps,x)| \leq C \ga_\eps^{|\al| - s}
\]
uniformly  in $x$. Hence another application of (\ref{Daub}) proves
the assertion.
\end{proof}

\newcommand{\SortNoop}[1]{}


\end{document}